\documentclass[12pt,draft]{article}

\begin{document}

\begin{center}
\LARGE\noindent\textbf{ Cycles of many lengths in digraphs with Meyniel-like condition }\\

\end{center}
\begin{center}
\noindent\textbf{Samvel Kh. Darbinyan}\\

Institute for Informatics and Automation Problems, Armenian National Academy of Sciences

E-mail: samdarbin@ipia.sci.am\\
\end{center}

\begin{center}
A detailed proof of a theorem published earlier by author in "Cycles of any length in digraphs with large semi-degrees", Academy Nauk Armyan SSR, Doklady 75(4) (1982) 147-152.   
\end{center}

\textbf{Abstract}

C. Thomassen (Proc. London Math. Soc. (3) 42 (1981), 231-251) gave a characterization of strongly connected non-Hamiltonian digraphs of order $p\geq 3$ with minimum degree $p-1$. In this paper we give an analogous characterization  of strongly connected non-Hamiltonian digraphs with Meyniel-type condition (the sum of degrees of every pair of non-adjacent vertices $x$ and $y$ at least $2p-2$). Moreover, we prove that such digraphs $D$ contain cycles of all lengths $k$, for $2\leq k\leq m$,
 where  $m$ is the length of a longest cycle in $D$.\\


\section{Introduction}

In this paper we consider finite digraphs (directed graphs)  without loops and multiple arcs. Every cycle and path is assumed simple and directed.

A  digraph $D$ is called {\it Hamiltonian} if it contains a Hamiltonian cycle , i.e., a cycle that includes all the vertices of $D$, and is {\it pancyclic} if it contains cycles of every length $m$, $3\leq m\leq p$, where $p$ is the number of the vertices of $D$.  
We recall the following well-known degree conditions (Theorems A and B), which guarantee that a digraph is Hamiltonian.\\

\noindent\textbf{Theorem A} (Ghouila-Houri \cite{[6]}). {\it Let $D$ be a strongly connected digraph of order $p\geq 2$. If $d(x)\geq p$ for all  vertices of $D$, then $D$ is Hamiltonian.}\\

For the next theorem we need the following definition.\\

\noindent\textbf{Definition.} {\it Let $D$ be a digraph of order $p$, and let  $k$ be an  integer. We will say that   $D$  satisfies condition $M_k$ when 
$$d(x)+d(y)\geq 2p-2+k$$ 
for all pairs of non-adjacent vertices $x$ and $y$ .}\\

\noindent\textbf{Theorem B} (Meyniel \cite{[9]}).
{\it Let $D$ be a strongly connected digraph of order $p\geq 2$ satisfying condition $M_1$. Then $D$ is Hamiltonian.}\\ 
 
Using Meyniel's theorem it is not difficult to show the following corollary.\\

\noindent\textbf{Corollary}.
{\it Let $D$ be a  digraph of order $p\geq 2$. If $d(x)+d(y)\geq 2p-3$ for every pair of non-adjacent vertices $x, y$ in $D$, 
then $D$ contains a Hamiltonian path}.\\

Theorems A and B are best possible.  Nash-Williams \cite{[10]} raised the following problem. \\

\noindent\textbf{Problem 1}: (Nash-Williams \cite{[10]}). {\it Describe all the extreme digraphs for Theorem A, i.e., describe all the
 strongly connected non-Hamiltonian digraphs of order $p$ with  minimum degree $p-1$.}\\ 

As a partial solution to the above problem of Nash-Williams, Thomassen  \cite{[13]} proved the following theorem.\\

\noindent\textbf{Theorem C} (Thomassen \cite{[13]}). {\it Let $D$ be a strongly connected non-Hamiltonian  digraph of order $p$  with minimum degree at least 
$p-1$. Let $S$ be a longest cycle in $D$. Then any two vertices of $V(D)\setminus V(S)$ are adjacent, every vertex of $V(D)\setminus V(S)$ has degree $p-1$ $($in $D$$)$, and every component of $D\langle V(D)\setminus V(S)\rangle$ is complete. Furthermore, if $D$ is 2-strongly connected, then $S$ can be chosen such that $D\langle V(D)\setminus V(S)\rangle$ is a transitive tournament.}\\

In \cite{[3]}, \cite{[4]}, \cite{[8]}, \cite{[11]} and \cite{[12]} are shown that: if a strongly connected digraph  $D$ satisfies the condition of Theorem A or the condition $M_i$, $1\leq i\leq 3$, then $D$ also is pancyclic unless $D$ is isomorphic to one of some exceptional digraphs which are characterized. In \cite{[5]}, the author and Mosesyan proved that if an orgraph $G$ (i.e., a digraph without cycles of length two) of order $2n+1\geq 17$ is $(n-1)$-biregular, then $G$
is pancyclic.  The reader can find more information on the Hamiltonian and pancyclic digraphs in the survey paper \cite{[1]} by Bermond and Thomassen.\\

Motivated of Meyniel's and Thomassen's theorems it is natural to consider the analogous problem for the Meyniel's theorem:\\

 \noindent\textbf{Problem 2}: {\it Characterize those strongly connected digraphs which satisfy condition $M_0$  but are not-Hamiltonian. Whether such digraphs $D$ contains cycles of  lengths $r$ for all $2\leq r\leq m$, where  $m$ is the length of the longest  cycles of $D$?}\\

We here prove a result analogous to the abovementioned theorem due to Thomassen for the Meyniel-like condition.  \\

\noindent\textbf{Theorem D}. {\it Let $D$ be a strongly connected non-Hamiltonian digraph of order $p\geq 3$ satisfying condition $M_0$. Let $C_m=x_1x_2\ldots x_mx_1$ be a longest cycle in $D$, and let $D_1$, $D_2, \ldots , D_h$ be the strong components of
$D\langle V(D)\setminus V(C_m)\rangle$ labelled in such a way that no vertex of $D_i$ dominates a vertex of $D_j$ whenever $i>j$.
Then the following statements hold:

I. Any two distinct vertices of $A:=V(D)\setminus V(C_m)$ are adjacent; every vertex of $A$ has degree at most $p-1$ in $D$; 
every component $D_i$ $($$1\leq i\leq h$$)$ is complete.

II. If $G\not\cong [(K_{p-m}\cup K_{m-1})+K_1]^*$, then for every $l\in [1,h]$ there are two distinct vertices $x_a,x_b$ on
 $C_m$  and  some vertices $u,v$ in $V(D_l)$ $($possibly $u=v$$)$ such that $x_au, vx_b\in E(D)$ and   
$$
E(B_l\rightarrow V(D_1)\cup V(D_2)\cup \ldots \cup V(D_l)) 
=E(V(D_l)\cup V(D_{l+1})\cup \ldots \cup V(D_h)\rightarrow B_l)=\emptyset,
$$
in particular, $E(V(D_l),B_l)=\emptyset$, where 
$B_l:=\{x_{a+1}, x_{a+2}, \ldots , x_{b-1}\}\not=\emptyset$. Moreover, $D\langle B_l\rangle$ is complete;  
$$
V(D_1)\cup V(D_{2})\cup \ldots \cup V(D_{l-1})\rightarrow  B_l\cup V(D_l)\rightarrow
V(D_{l+1})\cup \ldots \cup V(D_h);
$$
for all vertices $z\in V(D_l)$ and $y\in B_l$,  $d(z,C_m)=m-|B_l|+1$;  $d(y,C_m)=m+|B_l|-1$;  
$|B_l|\geq |V(D_l)|$ and
any vertex of $V(D_l)\cup B_l$ cannot be inserted into $C_m[x_b,x_a]$ $($in particular, 
$x_a\rightarrow B_l\cup V(D_l)\rightarrow x_b$$)$.

III. If $D$ is 2-strongly connected, then the induced subdigraph $D\langle A\rangle$ is a transitive tournament.}

IV. $D$ contains  cycles of every length $r$,  $r\in [2,m]$, unless when $p$  odd  and $D$ is isomorphic to the complete bipartite digraph $ K^*_{\lfloor p/2\rfloor,\lfloor p/2\rfloor+1}$.\\ 

It is worth remaking that in the proofs of the first and second statements of Theorem D we use some ideas appeared in \cite{[13]}.
The aim of this paper is to present a detailed proof of Theorem D.

\section{Terminology and notation}

  In this paper we consider finite digraphs without loops and multiple arcs. 
We refer the reader to \cite{[1]} and \cite{[7]}.
for terminology not discussed here. For a digraph $D$, we denote by $V(D)$ the vertex set of $D$ and by $E(D)$ the set of arcs in $D$. The {\it order} of $D$ is the number of its vertices. 

  The arc of a digraph $D$ directed from
   $x$ to $y$ is denoted by $xy$ (we say that {\it $x$ dominates $y$}).  For disjoint subsets $A$ and  $B$ of $V(D)$  we define $E(A\rightarrow B)$ \,
   as the set $\{xy\in E(D) | x\in A, y\in B\}$ and $E(A,B)=A(A\rightarrow B)\cup A(B\rightarrow A)$. If $x\in V(D)$
   and $A=\{x\}$ we write $x$ instead of $\{x\}$. If $A$ and $B$ are two disjoint subsets of $V(D)$ such that every
   vertex of $A$ dominates every vertex of $B$, then we say that {\it $A$ dominates $B$}, denoted by $A\rightarrow B$. 
The {\it out-neighborhood} of a vertex $x$ is the set 
$O(x)=\{y\in V(D) | xy\in E(D)\}$ and $I(x)=\{y\in V(D) | yx\in E(D)\}$ is the {\it in-neighborhood} of $x$. Similarly, if $A\subseteq V(D)$, then $O(x,A)=\{y\in A | xy\in E(D)\}$ and 
$I(x,A)=\{y\in A | yx\in E(D)\}$. 
 The {\it degree} of the vertex $x$ in $D$ is defined as $d(x)=|O(x)|+|I(x)|$ (similarly,
 $d(x,A)=|O(x,A)|+|I(x,A)|$).

 For integers $a$ and $b$, $a\leq b$, let $[a,b]$  denote the set of
all integers which are not less than $a$ and are not greater than
$b$. We denote $left[a,b]=a$ and $right[a,b]=b$.

 The path (respectively, the cycle) consisting of the distinct vertices $x_1,x_2,\ldots ,x_m$ ($m\geq 2 $) and the arcs $x_ix_{i+1}$, $i\in [1,m-1]$  (respectively, $x_ix_{i+1}$, $i\in [1,m-1]$, and $x_mx_1$), is denoted by  $x_1x_2\cdots x_m$ (respectively, $x_1x_2\cdots x_mx_1$). 
We say that $x_1x_2\cdots x_m$ is a path from $x_1$ to $x_m$ or is an {\it $(x_1,x_m)$-path}. The {\it length} of a cycle (of a path) is the number of its arcs. $C_k$ ($k\geq 2$) will denote the cycle of length $k$. A cycle (respectively, a path) that contains  all the vertices of $D$  is a {\it Hamiltonian cycle} (respectively, a {\it Hamiltonian path}). A digraph $D$ is {\it Hamiltonian} if it contains a Hamiltonian cycle.
 
If $P$ is a path containing a subpath from $x$ to $y$, then  $P[x,y]$ denotes the subpath of $P$ from $x$ to $y$.
 Similarly, if $C$ is a cycle containing vertices $x$ and $y$, then $C[x,y]$ denotes the subpath of $C$ from $x$ to $y$. For convenience, we also use $P[x,y]$ ($C[x,y]$) to denote the vertex set of the corresponding subpath.

Given a vertex $x$ of a directed path $P$ or  a directed cycle $C$, we use the notation $x^+$ and $x^-$ for the successor and the predecessor of $x$ (on $P$ or on $C$) according to the orientation. 

The subdigraph of $D$ induced by a subset $A$ of $V(D)$ is denoted by $D\langle A\rangle$. A digraph $D$ is {\it strongly connected} (or, just, {\it strong}) if there exists a path from $x$ to $y$ and a path from $y$ to $x$ for every pair of distinct vertices $x,y$. A digraph $D$ is {\it $k$-strongly connected} (or {\it $k$-strong}, $k\geq 1$) if $|V(D)|\geq k+1$ and $D\langle V(D)\setminus A\rangle$ is strong for any set $A$ of at most $k-1$ vertices. 
By Menger's theorem, this is equivalent to the property that for any ordered pair of distinct vertices $x,y$ there are $k$ internally disjoint paths from $x$ to $y$. 

 A {\it strong component}  of a digraph $D$ is a maximal induced strong subdigraph of $D$. Two distinct vertices $x$ and $y$ in $D$ are adjacent if $xy\in A(D)$ or $yx\in A(D) $ (or both). 
  For an undirected graph $G$, we denote by $G^*$ the symmetric digraph obtained from $G$ by replacing every edge $xy$ with the pair $xy$, $yx$ of arcs.
  
  We will denote the complete bipartite digraph with partite sets of cardinalities $p$, $q$ by $K^*_{p,q}$.

\section{Preliminaries}

Let us recall some well-known lemmas used in this paper.\\

\noindent\textbf{Lemma 1} (\cite{[8]}). {\it Let $D$ be a digraph of order $p\geq 3$  containing a cycle $C_m$, $m\in [2,p-1]$. Let $x$ be a vertex not contained in this cycle. If $d(x,V(C_m))\geq m+1$, then for every $k$, $k\in [2,m+1]$, $D$ contains a cycle $C_k$ including $x$}. 
 \\

The following lemma  is a  modification of a lemma by Bondy and Thomassen  \cite{[2]}.\\

\noindent\textbf{Lemma 2}. {\it Let $D$ be a digraph of order $p\geq 3$  containing a path $P:=x_1x_2\ldots x_m$, $m\in [2,p-1]$.    Let $x$ be a vertex not contained in this path. If one of the following statements holds:

 (i) $d(x,V(P))\geq m+2$; 

 (ii)  $d(x,V(P))\geq m+1$ and $xx_1\notin E(D)$ or $x_mx\notin E(D)$; 

 (iii)  $d(x,V(P))\geq m$, $xx_1\notin E(D)$ and $x_mx\notin E(D)$};

\noindent\textbf{}{\it then there is an  $i\in [1,m-1]$ such that $x_ix,xx_{i+1}\in E(D)$, i.e., $D$ contains a path $x_1x_2\ldots x_ixx_{i+1}\ldots x_m$ of length $m$}  ({\it we say that  $x$ can be inserted into $P$ or the path  $x_1x_2\ldots x_ixx_{i+1}\ldots x_m$ is an extended path obtained from $P$ with  $x$}). \\

As a consequence of Lemma 2, we get the following Lemma 3 (we give the proof of Lemma 3 here for completeness).\\

\noindent\textbf{Lemma 3}. {\it Let $D$ be a digraph of order $p\geq 4$. Suppose that    $P:=x_1x_2\ldots x_m$, $m\in [2,p-2]$, is a longest path from $x_1$ to $x_m$ in $D$. If the induced subdigraph $D\langle V(D)\setminus V(P)\rangle$ is strongly connected and   for all vertices $x$ of $ V(D)\setminus V(P)$, $d(x,V(P))=m+1$,  then there is an integer $l\in [1, m]$ such that
 $$
O(x,V(P))=\{x_1,x_2,\ldots ,x_l\} \quad \hbox{and} \quad I(x,V(P))=\{x_l,x_{l+1},\ldots , x_m\}.
$$}  

\noindent\textbf{Proof}. Let $x$ be an arbitrary vertex of $A:=V(D)\setminus V(P)$. Since $d(x,V(P))=m+1$ and $x$ cannot be inserted into $P$, by Lemma 2(ii) we have $xx_1$ and $x_mx\in E(D)$. We claim that $x$ is adjacent to every vertex of $P$. Assume that this is not the case. Let $x$ and $x_k$, where $2\leq k\leq m-1$, are not adjacent. Put $P_1:=x_1x_2\ldots x_{k-1}$ and $P_2:=x_{k+1}x_{k+2}\ldots x_m$ (possibly, $k=2$ or $k=m-1$). Since $x$ cannot be inserted neither into $P_1$ nor in $P_2$, by Lemma 2(i) we have that $d(x,V(P_1))\leq |V(P_1)|+1$ and $d(x,V(P_2))\leq |V(P_2)|+1$. On the other hand,
$$
m+1=d(x,V(P))=d(x,V(P_1))+d(x,V(P_2))\leq k+m-k+1=m+1.
$$
This implies that $d(x,V(P_1))=|V(P_1)|+1$ and $d(x,V(P_2))= |V(P_2)|+1$. Again applying Lemma 2(ii) to $P_1$ and to $P_2$, we obtain $x_{k-1}x$ and $xx_{k+1}\in E(D)$. Since $D\langle A\rangle$ is strong and $P$ is a longest $(x_1,x_m)$-path, it follows that $x_k$ and every vertex $z$ of $A\setminus \{x\}$ are not adjacent. Now using  Lemma 2, by similar arguments, we conclude that $x_{k-1}z$ and $zx_{k+1}\in E(D)$ since $d(z,V(P))= m+1$. Hence, it is not difficult to describe an $(x_1,x_m)$-path of length greater than $m-1$, which is a contradiction. 

Thus we have proved that every vertex of $A$ is adjacent to  every vertex of $V(P)$. Then, since $D\langle A\rangle$ is strong and  $d(x,V(P))=m+1$ for all $x\in A$, there exists an integer $l\in [1,m]$ such that 
$O(x,V(P))=\{x_1,x_2,\ldots ,x_l\}$ and $I(x,V(P))=\{x_l,x_{l+1},\ldots , x_m\}$. Lemma 3 is proved. \fbox \\\\

\section{Proof of Theorem D}

 For any integers $i$ and $k$ ($1\leq i\leq k\leq h$) put $A_i:=V(D_i)$, $A_{i,k}:=\cup ^k_{j=i}A_j$, $a_i:=|A_i|$ and $a_{i,k}:=|A_{i,k}|$.  We use $C_m$  also for $V(C_m)$.  Since $C_m$ is a longest cycle in $D$, using Lemma 1, we obtain that for every vertex $y$ of $A$,  $d(y,C_m)\leq m$. This together with  condition $M_0$ implies that for any two  non-adjacent distinct vertices $y$ and $z$ of $A$ the following holds
$$
2p-2\leq d(y)+d(z)=d(y,C_m)+d(z,C_m)+d(y,A)+d(z,A)\leq 2m+d(y,A)+d(z,A).
$$
Hence, $d(y,A)+d(z,A)\geq 2(p-m)-2$, i.e., the subdigraph  $D\langle A\rangle$ satisfies condition $M_0$. Therefore, by corollary of Meyniel's theorem, the subdigraph $D\langle A\rangle$  has a Hamiltonian path. In particular, for each $i\in [1,h-1]$ there is an arc from a vertex of $A_i$ to a vertex of $A_{i+1}$, i.e., 
$E(A_i\rightarrow A_{i+1})\not= \emptyset$. 
From this and strongly connectedness of $D$ it follows that  
$$
E(C_m\rightarrow A_{1})\not= \emptyset \quad \hbox{and} \quad  E(A_{h}\rightarrow C_m)\not= \emptyset. \eqno (1)
$$
We consider two cases.

\noindent\textbf{Case 1}. Exactly  one vertex, say $x$, of the cycle $C_m$ is adjacent to some vertices of $ A$.

Then $E(A, C_m\setminus \{x\})=\emptyset$. By condition $M_0$, for every pair of vertices $y\in C_m\setminus \{x\}$ and $z\in A$ we
have 
$$
2p-2\leq d(y)+d(z)=d(y,C_m)+d(z,A)+d(z,\{x\}).  \eqno (2)
$$
Since
 $$d(z,\{x\})\leq 2, \quad d(y,C_m)\leq 2m-2 \quad \hbox{and} \quad d(z,A)\leq 2(p-m-1),  \eqno (3)
$$
from (2) it follows that all the inequalities of (3), in fact, are equalities. Therefore, the subdigraphs  $D\langle C_m\rangle$ and $D\langle A\cup \{x\}\rangle$ are complete digraphs. This means that  $D\cong [(K_{p-m}\cup K_{m-1})+K_1]^*$. In particular, $D$ is not 2-strongly connected. It is easy to see that $m\geq p-m+1$ since $C_m$ is a longest cycle in $D$. From this and $d(z)=2p-2m$, where $z\in A$, we have $d(z)=2p-2m\leq p-1$. Thus, we have proved that in this case the theorem is true.\\

\noindent\textbf{Case 2}. There are at least two distinct vertices  on $C_m$ which are  adjacent to some vertices of $A$.\\

We first prove   Claim 1 and Claim 2.\\

\noindent\textbf{Claim 1}. {\it Let $1\leq l\leq q\leq h$. Suppose that there are  vertices $x_a$, $x_b$ with $x_a\not=x_b$  on $C_m$, and vertices 
$u\in A_l$, $v\in A_q$ $($possibly,  $u=v$$)$ such that $x_au, vx_b\in E(D)$.   
Moreover, assume that for some integer $k\in [l,q]$ and $B_k:=\{x_{a+1},x_{a+2},\ldots , x_{b-1}\}$,
$$
E(B_k\rightarrow A_{1,k})= E(A_{k,h}\rightarrow B_k)= \emptyset. \eqno (4)
$$
 Then the following statements hold:

a. $|B_k|\geq 1$, $l=k=q$, the subdigraphs  $ D_k$ and $D\langle B_k\rangle$ are complete digraphs;

b. For all vertices $y\in B_k$ and $z\in A_k$, $d(y,C_m)=m+|B_k|-1$, $d(z,C_m)=m-|B_k|+1$ and $d(z)\leq p-1$;  

c. The path $C_m[x_b,x_a]$ cannot be extended with any vertex of $A_k\cup B_k$ $($in particular, $x_a\rightarrow A_k\cup B_k\rightarrow x_b$$)$  and  
$A_{1,k-1}\rightarrow A_k\cup B_k\rightarrow A_{k+1,h}$}.

\noindent\textbf{Proof of Claim 1}. From (4) it follows that $E(A_k,B_k)=\emptyset$. Since $C_m$ is a longest cycle in $D$ and since in $D\langle A\rangle$ there is an $(u,v)$-path,  whose existence (if $l\not= q$) follows from the fact that 
$D\langle A\rangle $ has a Hamiltonian path,
 we have $|B_k|\geq 1$. Now we extend the path $C_m[x_b,x_a]$ with 
the vertices of $B_k$ as much as possible. We obtain an $(x_b,x_a)$-extended path, say $Q$. Because of the maximality of $C_m$, the presence of the arcs $x_au$, $vx_b$ and an $(u,v)$-path in $D\langle A\rangle $,  some vertices $y_1,y_2,\ldots , y_d$ of $B_k$, where $1\leq d\leq |B_k|$, are not on the obtained extended path $Q$. Now using Lemma 2, we obtain
$$
d(y_i,C_m)=d(y_i,V(Q))+d(y_i,\{y_1,y_2,\ldots , y_d\})\leq m+d-1. \eqno (5)
$$
 
Let $z$ be an arbitrary vertex of $A_k$. It is easy to see that the vertex $z$ cannot be inserted into $C[x_b,x_a]$. 
Now using the fact that $E(A_k,B_k)=\emptyset$ and Lemma 2, we obtain
$$
d(z,C_m)\leq m-|B_k|+1. \eqno (6)
$$
On the other hand, by (4), we have
$$
d(y_i)+d(z)=d(y_i,C_m)+d(z,C_m)+|E(A_{1,k-1}\rightarrow y_i)|+|E(y_i\rightarrow A_{k+1,h})| $$ $$+|E(A_{1,k}\setminus \{z\}\rightarrow z)|+|E(z\rightarrow A_{k,h}\setminus \{z\})|.
$$
Now using (5), (6) and condition $M_0$, we get
$$
2p-2\leq d(y_i)+d(z)\leq m+d-1+m-|B_k|+1+a_{1,k-1}+a_{k+1,h}
$$ 
$$
+a_{1,k}+a_{k,h}-2=2p+d-|B_k|-2                    \eqno (7)
$$
because of $a_{1,k}+a_{k+1,h}=a_{1,k-1}+a_{k,h}=p-m$. Since $d\leq |B_k|$, we have that $d= |B_k|$, i.e., the path 
$C_m[x_b,x_a]$ cannot be extended with any vertex of $B_k$. Moreover, it follows that there must be equalities in all estimates that led to (7) as well, i.e., the subdigraphs $D_k$ and $D\langle B_k\rangle $ are complete digraphs, 
$d(z,C_m)=m-|B_k|+1$ and  $d(y,C_m)=m+|B_k|-1$ and 
$$
A_{1,k-1}\rightarrow A_k\cup B_k\rightarrow A_{k+1,h}.
$$
From the last expression  it follows that $u,v\in A_k$, i.e., $l=q=k$. Since $d(z,C_m)=m-|B_k|+1$ and $d(z,B_k)=0$, using Lemma 2
we obtain that $x_a\rightarrow A_k\rightarrow x_b$. Similarly, $x_a\rightarrow B_k\rightarrow x_b$. From this, the arbitrariness of $z$ and the fact that $D_k$ is complete it follows that $|B_k|\geq a_k$, which in turn implies that 
$$
d(z)=d(z,C_m)+d(z,A_k)+d(z,A\setminus A_k)= m-|B_k|+1+2a_k-2+p-m-a_k\leq p-1.
$$
This completes the proof of Claim 1. \fbox \\\\

\noindent\textbf{Claim 2}. {\it Suppose that $l\in [1,h]$. Then there exist two distinct vertices $x_a$, $x_b$ on $C_m$ and  some vertices 
$u, v$ in $A_l$  $($possibly $u=v$$)$ such that $x_au, vx_b\in E(D)$ and     
$$
E(B_l\rightarrow A_{1,l})= E(A_{l,h}\rightarrow B_l)= \emptyset, 
$$
where  $B_l:=\{x_{a+1},x_{a+2},\ldots , x_{b-1}\}$ and $|B_l|\geq 1$ $($i.e., for every $l\in [1,h]$ the suppositions of Claim 1 hold when $l=q=k$ $)$}.

\noindent\textbf{Proof of Claim 2}. In order to prove Claim 2, we first prove that Claim 2 is true for $l=1$ or $l=h$.

 Since $C_m$ is a longest cycle in $D$ and since $D\langle A\rangle$ contains a Hamiltonian path, from (1) it follows that there exists a vertex $x_i\in C_m$ such that  $E(x_i\rightarrow A_1)\not=\emptyset$, 
$E(x_{i+1}\rightarrow A_1)=\emptyset$ and  $E(A\rightarrow x_{i+1})=\emptyset$.

 From  $E(A_h\rightarrow C_m)\not =\emptyset$ (by (1))
 and
$E(A\rightarrow x_{i+1})=\emptyset$ it follows that there is a vertex $x_b\in C_m\setminus \{x_{i+1}\}$ such that
 $$
E(A\rightarrow x_b)\not =\emptyset \quad \hbox{and} \quad E(A\rightarrow \{x_{i+1},\ldots ,  x_{b-1}\}) =\emptyset. \eqno (8) 
$$
To be definite, assume that $zx_b\in E(D)$, where $z\in A_q$ and $q\in [1,h]$. 

Assume first that the vertices $x_i$ and $x_b$ are distinct. It is easy to see that $E(x_{b-1}\rightarrow A_{1,q}) =\emptyset$ since $C_m$ is a longest cycle in $D$. This together with (8) implies that
 \\ $E(x_{b-1},A_{1,q}) =\emptyset$. Now from 
$E(x_i\rightarrow A_1)\not=\emptyset$ it follows that there is a vertex $x_a\in \{x_i,x_{i+1},\ldots , x_{b-2}\}$ and an integer  $k\in [1,q]$ such that 
$$
E(x_a\rightarrow A_k)\not =\emptyset \quad \hbox{and} \quad E(\{x_{a+1},\ldots ,  x_{b-1}\}\rightarrow A_{1,q}) =\emptyset. \eqno (9) 
$$
By (8) and (9), for $B_l:= \{x_{a+1},\ldots ,  x_{b-1}\}$ we have that
$$
 E(B_l\rightarrow A_{1,q}) = E(A\rightarrow B_l)=\emptyset,  
$$
in particular, $E(A_{k,q}, B_l)=\emptyset$. Thus, the suppositions of Claim 1 hold, which in turn implies that $k=q$ and 
$A_{1,k-1}\rightarrow B_l\cup A_k\rightarrow A_{k+1,h}$. 
Therefore, $A_{1,k-1}=\emptyset$, i.e., $k=1$, since $E(A\rightarrow B_l)=\emptyset$ by (8).
Thus in this case ($x_i\not=x_b$) for $l=1$ Claim 2 is true. 

Assume next that  there is no $x_j$ other than $x_i$ such that 
$E(A\rightarrow x_j)\not=\emptyset$. We have  $E(A_h\rightarrow C_m)\not=\emptyset$ (by (1)) and $E(A\rightarrow C_m\setminus \{x_i\})=\emptyset$, in particular, $E(A_h\rightarrow x_i)\not=\emptyset$. Under the condition of Case 2,  there exists a vertex $x_g\in C_m$ other than $x_i$ such that 
$E(x_g\rightarrow A)\not=\emptyset$. 
It is not difficult to check that in the converse digraph of $D$ the considered case $x_i\not= x_b$ holds. In this case Claim 2 is true for $l=h$. So, Claim 2 is true for $l=1$ or $l=h$. This means that if $h=1$, then Claim 2 is proved. Assume that $h\geq 2$. 

Without loss of generality we may assume that Claim 2 is true for $l=h$ 
 (if Claim 2 is true for $l=1$, then we will consider the converse digraph of $D$).

Now we assume that Claim 2 is true for $t+1$, where $2\leq t+1\leq h$, and prove it for $t$. Then  there are vertices $u',v'\in A_{t+1}$ and two distinct vertices $x_d,x_s$ on $C_m$ such that $x_du', v'x_s\in E(D)$ and 
 $$
E(\{ x_{d+1},x_{d+2},\ldots , x_{s-1}\}\rightarrow A_{1,t+1})= E(A_{t+1,h}\rightarrow \{x_{d+1},x_{d+2},\ldots , x_{s-1}\})= \emptyset.
$$ 
Then, by Claim 1, $A_{1,t}\rightarrow 
\{x_{d+1},x_{d+2},\ldots , x_{s-1}\}\cup A_{t+1} \rightarrow A_{t+2,h}  
$
in particular, $A_{1,t}\rightarrow \{x_{d+1}\}$ which in turn implies that $A_{t}\rightarrow \{x_{d+1}\}$. This implies that $E(x_d\rightarrow A_{1,t})=\emptyset$, since $C_m$ is a longest cycle in $D$. 
This together with $E(C_m\rightarrow A_1)\not=\emptyset$ implies that there exists a vertex
 $x_i\notin \{x_d,x_{d+1}, \ldots , x_{s-1}\}$ such that 
$$
E(x_i\rightarrow A_{1,t})\not=\emptyset \quad \hbox{and} \quad E(\{x_{i+1},x_{i+2},\ldots , x_{d}\} \rightarrow A_{1,t})=\emptyset.         \eqno (10)
$$

To be definite, assume that $x_iy\in E(D)$, where $y\in A_r$ and $r\in [1,t]$. Since $C_m$ is a longest cycle in $D$,  it follows that $E(A_{r,h}\rightarrow x_{i+1})=\emptyset $, in particular 
 $E(x_{i+1},A_{r,t})=\emptyset$ by (10). This together with $A_t\rightarrow x_{d+1}$  implies that  there exists a  vertex  
$x_j\in \{x_{i+2},x_{i+3},\ldots , x_{d+1}\}$ such that 
$$
E(A_{r,h}\rightarrow  x_{j})\not=\emptyset\quad \hbox{and} \quad 
E(A_{r,h}\rightarrow \{x_{i+1},x_{i+2},\ldots , x_{j-1}\})=\emptyset. \eqno (11)
$$
To be definite, assume that $E(A_{q}\rightarrow  x_{j})\not=\emptyset$, where $r\leq q\leq h$.

Assume that $q\geq t$, then by (10) and (11) we have
$$
E(\{x_{i+1},x_{i+2},\ldots , x_{j-1}\}\rightarrow A_{1,t})= E(A_{r,h}\rightarrow \{x_{i+1},x_{i+2},\ldots , x_{j-1}\}=\emptyset.
$$
Therefore, if $q\geq t$, then from  $r\leq t$ and Claim 1 it follows that $l=t=q$, which in turn implies that for $t$  Claim 2 is true.

We may therefore assume  that $q\leq t-1$. Let $B_k:=\{x_{i+1},x_{i+2},\ldots , x_{j-1}\}$ and let $z$ be an arbitrary vertex of $A_k$, 
where $k\in [r,q]$. (10) and (11), in particular, mean that $E(z,B_k)=\emptyset$. Since $C_m$ is a longest cycle in $D$, the vertex $z$ cannot be inserted into $C_m$. Using Lemma 2, 
 and $E(z,B_k)=\emptyset$ 
 we obtain that $d(z, C_m)\leq m-|B_k|+1$. It is clear that $d(z,A)\leq p-m+a_k-2$. Thus, 
$d(z)\leq p-|B_k|+a_k-1$.

Since $x_iy, wx_j\in E(D)$, where $y\in A_r$, $w\in A_q$ and $1\leq r\leq q\leq t-1$, and since $D\langle A\rangle$ has a Hamiltonian 
path (in particular, in $D\langle A\rangle$ there is a $(y,w)$-path), it follows that the path $C_m[x_j,x_i]$ cannot be extended with all  vertices of $B_k$. This means that for some vertices $y_1,y_2,\ldots , y_d$ of $B_k$ ($d\geq 1$) the following holds $d(y_i,C_m)\leq m+d-1$. Using (10) and (11), we obtain that $E(y_i\rightarrow A_{1,t})=E(A_{r,h}\rightarrow y_i)=\emptyset$. Therefore, $A(y_i, A_{r,t})=\emptyset$ and
$$
d(y_i,A)= d(y_i,A_{1,r-1})+d(y_i,A_{t+1,h})\leq a_{1,r-1}+a_{t+1,h}
$$
(if $r=1$, then $A_{1,r-1}=\emptyset)$.
 Now, since the vertices $z$ and $y_i$ are not adjacent, condition $M_0$ implies that 
$$
2p-2\leq d(y_i)+d(z)\leq p-|B_k|+a_k-1+a_{1,r-1}+a_{r,t}+a_{t+1,h}-a_{r,t} $$ $$ + m+d-1\leq  2p-|B_k|+d-2+a_k-a_{r,t},
$$
since $p-m=a_{1,r-1}+a_{r,t}+a_{r+1,h}$. Using the facts that $a_k\leq a_{r,t}$ and $d\leq |B_k|$, we obtain that $d=|B_k|$ and
$a_k=a_{r,t}$.
 The last equality is possible if $k=r=t=q$, which contradicts the assumption that $q\leq t-1$. This completes the proof of Claim 2.  \fbox \\\\

 The first and second statements of the theorem in Case 2 follows immediately from Claims 1 and 2.\\

Now we will prove the third statement of the theorem. 
Suppose in addition that $D$ is 2-strongly connected. We want to prove that the induced subdigraph $D\langle A\rangle$ is a transitive tournament. By the first statement of the theorem, it suffices to prove that $a_i=1$ for all $i\in [1,h]$. 
Assume that this is not the case. Then for some $k\in [1,h]$, $a_k\geq 2$. By the first and second statements of the theorem, the subdigraphs $D_k$ and 
$D\langle B_k\rangle$  are complete digraphs, where $B_k=\{x_{a+1},\ldots , x_{b-1}\}$ and $x_a\rightarrow D_k\rightarrow x_b$. Therefore, $s:=b_k\geq a_k\geq 2$ (recall that $b_k=|B_k|$ and $a_k=|A_k|$). Put $P:=C_m[x_b, x_a]=x_1x_2\ldots x_{m-s}$. Let $y$ (respectively, $z$) be an arbitrary vertex of $B_k$ (respectively, $A_k$). By the second statement of the theorem,  
any vertex of $A_k\cup B_k$ cannot be inserted into $P$ and $d(y,V(P))=d(z,V(P))=m-s+1$.  Therefore, by Lemma 3, there exist $l, l'\in [1,m-s]$ such that
$$
O(z,V(P))=\{x_1,x_2,\ldots  ,x_{l}\}, \quad I(z,V(P))=\{x_{l},x_{l+1},\ldots , x_{m-s}\} \eqno (12) $$
$$
O(y,V(P))=\{x_1,x_2,\ldots  ,
x_{l'}\}, \quad I(y,V(P))=\{x_{l'},x_{l'+1},\ldots , x_{m-s}\}. \eqno (13)
$$
Without loss of generality, we may assume that $l\leq l'$ (for otherwise we will consider the converse digraph of $D$).
Since $C_m$ is a longest cycle in $D$, it is easy to see that
$$
E(\{x_1,x_2,\ldots , x_{l-1}\}\rightarrow A_{1,k})=E(A_{k,h}\rightarrow \{x_{l+1}, x_{l+2},\ldots , x_{m-s}\})
=\emptyset. \eqno (14)
$$
 If $x_ix_j\in E(D$ with $i\in [1,l-1]$ and $j\in [l+1,m-s]$, then, by (12) and (13),  the cycle $x_1x_2\ldots x_ix_j\ldots  x_{m}$ $x_{i+1}\ldots$ $ x_{j-1}zx_1$, where $z\in A_k$, is a cycle of length $m+1$, which contradicts that $C_m$ is a longest cycle in $D$.  We may therefore assume that
$$
E(\{x_1,x_2,\ldots , x_{l-1}\}\rightarrow \{x_{l+1},,x_{l+2},\ldots , x_{m-s}\})
=\emptyset.        \eqno (15)
$$
By the second statement of the theorem we have 
$$
E(B_k\rightarrow A_{1,k})=E(A_{k,h}\rightarrow B_k)=\emptyset.  \eqno (16)
$$
 This together with (12)-(15) implies  that if $2\leq l\leq m-s$, then 
$$
E(\{x_1,x_2,\ldots , x_{l-1}\}\cup A_{k+1,h}\rightarrow B_k\cup A_{1,k}\cup \{x_{l+1},\ldots , x_{m-s}\})=\emptyset.
$$
This means that $D-\{x_{l}\}$ is not strong,  which contradicts the assumption that $D$ is 2-strongly connected. 

Thus, we may assume  that $l=1$.  
   Let $l'\geq 2$, then using (13) it is not difficult to see that 
$$
E(x_1\rightarrow \{x_{3},x_{4},\ldots , x_{m-s}\}) 
=E(\{x_1\rightarrow A_{1,k-1})=\emptyset.       \eqno (17)
$$
Indeed, if $x_1x_i\in E(D)$ with $i\in [3,m-s]$, then the 
cycle $x_1x_i\ldots x_{m-s}\ldots x_mx_2\ldots x_{i-1}zx_1$  is a cycle of length $m+1$, and if $x_1u\in E(D)$, where 
$u\in A_{1,k-1}$, then, by the second statement of the theorem,  $A_{1,k-1}\rightarrow B_k$ (if $k\geq 2$) and the cycle 
$x_1ux_{m-s+1}\ldots x_mx_2\ldots x_{m-s}zx_1$ is a cycle of length $m+1$, a contradiction.

 Now using (14), (16) and (17), we obtain
$$
E(\{x_1\}\cup A_{k,h}\rightarrow B_k\cup A_{1,k-1}\cup \{x_{3},x_4, \ldots , x_{m-s}\})=\emptyset.
$$
This means that
 the subdigraph
$D-\{x_2\}$ is not 2-strongly connected, which is a contradiction.

Let finally $l=l'=1$. Then we have
$E(A_{k,h}\rightarrow B_k\cup A_{1,k-1}\cup \{x_{2},x_{3},\ldots , x_{m-s}\})
=\emptyset$. 
Therefore, $G-\{x_1\}$ is not strongly connected, which contradicts that $D$ is 2-strongly connected. This completes the proof of the third statement of the
theorem. \\

Finally we will prove the fourth statement of the theorem
 that for every $r\in [2,m]$, $D$ contains a cycle of length $r$
unless $p$ odd and $D$ isomorphic to the  complete bipartite digraph 
$K^*_{\lfloor {p/2}\rfloor , \lfloor {p/2}\rfloor+1}$.

Assume  first that there exists an integer $k\in [1,h]$ for which $|B_k| \geq 2$ (Claim 2). Put $P_k:=C_m[x_b,x_a]:=x_1x_2\ldots x_{m-s}$, where $s=|B_k|$. By Claims 1 and 2 we have, $D\langle B_k\rangle$ is a complete digraph, any vertex $y$ of $B_k$
cannot be inserted into $P_k$ and $d(y,V(P_k))=m-s+1$. Therefore, by Lemma 3, there exists an integer $l\in [1,m-s]$ 
such that
$$
\{x_{l}, x_{l+1},\ldots , x_{m-s}\}\rightarrow B_k\rightarrow \{x_1, x_{2},\ldots , x_{l}\}.
$$
Hence it is not difficult to check that $D$ contains a cycle of length $r$ for every $r\in [2,m]$.

In what follows we assume that $|B_i|=1$ for all $i\in [1,h]$. From the second statement  of the theorem it follows that $a_i=1$  for every $i\in [1,h]$. 
 Then, by the first statement of the theorem, 
in $D\langle A\rangle$ any two vertices are adjacent, which in turn implies  that  $D\langle A\rangle$ is a transitive tournament. 

Put $A_1:=\{x\}$. By Claim 1, we have that $d(x,C_m)=m-|B_1|+1$, i.e., $d(x,C_m)=m$ since $|B_1|=1$. Now for the cycle $C_m$ and the vertex $x$ it is not difficult to show the following proposition.\\

\noindent\textbf{Proposition 1}. Let $x_i$ be an arbitrary vertex on $C_m$. Then the following holds:

(i) {\it if $x_ix\notin E(D)$, then $xx_i^+\in E(D)$;}

(ii) {\it if $xx_i\notin E(D)$, then $x_i^-x\in E(D)$;}

(iii) {\it if $E( x,x_i)=\emptyset$, then $x_i^-x$ and $xx_i^+\in E(D)$.}

\noindent\textbf{Proof of Proposition 1}. Indeed, if $x_ix\notin E(D)$ and $xx_i^+\notin E(D)$, then, since $d(x,C_m)=m$ and  $C_m$ is a longest cycle in $D$, using Lemma 2(iii), we obtain 
$$
m=d(x,C_m)=d(x, C_m[x_i^+,x_i])\leq  |C_m[x_i^+,x_i]|-1=m-1,
$$
a contradiction.
In the same way, one can show that (ii) also is true. (iii) is an immediate consequence of  (i) and (ii). \fbox \\\\

   Let $x_{t_1},x_{t_2},\ldots , x_{t_n}$ be the vertices of $C_m$ that are not adjacent to $x$  numbered along the orientation of the cycle $C_m$. By Proposition 1(iii) we have
$$
x_{t_i}^-x\in E(D) \quad \hbox{and} \quad xx_{t_i}^+\in E(D).                   \eqno (18)
$$ 
Observe that any path $Q_i:=x_{t_i+1}x_{t_i+2}\ldots x_{t_{i+1}-1}$ (in other words, $Q_i=C_m[x^+_{t_i},x^-_{t_{i+1}}]$) cannot be extended with the vertex $x$
 (here, $t_{n+i}=t_i$ for all $i\in [1,n]$). Then, since the vertex $x$ is adjacent to any vertex of $V(Q_i)$, using Lemma 2
and the fact that $d(x,C_m)=m$, we obtain that $d(x,V(Q_i))=|V(Q_i)|+1$ and there exists an integer $l_i\geq 1$ such that
$$
O(x,V(Q_i))=\{x_{t_i+1},x_{t_i+2},\ldots , x_{t_i+l_i}\}, \quad I(x,V(Q_i))=\{x_{t_i+l_i},x_{t_i+l_i+1},\ldots , x_{t_{i+1}-1}\}. \eqno (19) 
$$
From (18) and the second statement of the theorem  it follows that
$$
d(x_{t_i},C_m)=m, \quad d(x_{t_i})=p-1, \quad x_{t_i}\rightarrow A\setminus \{x\} \quad \hbox{and} \quad E( A\rightarrow x_{t_i})=\emptyset.           \eqno (20)
$$
From (19) it follows that if $n=1$, then for every $r\in [2,m]$, $D$ contains a cycle of length $r$. In the sequel, we assume that $n\geq 2$, i.e., the number of vertices on $C_m$ which are not adjacent to $x$ is more than or equal to two.\\   

We need to prove the following Claims 3-6.\\

\noindent\textbf{Claim 3}. {\it Suppose that two distinct vertices $x_a$,  $x_b$ of 
$\{x_{t_1},x_{t_2},\ldots , x_{t_n}\}$ are not adjacent, then  the arcs
$x_{b-1}x_a$, $x_ax_{b+1}$, $ x_bx_{a+1}$,  $x_{a-1}x_b$ are in $E(D)$}.

\noindent\textbf {Proof of Claim 3}. By the first equality of (20), we have $d(x_a,C_m)=d(x_b,C_m)=m$. By (18), $x_{b-1}x$ and $xx_{b+1}\in E(D)$. This implies that   $x_b$ cannot be inserted  into $C_m[x_{b+1},x_{b-1}]$ since $C_m$ is a longest cycle in $D$. Hence, using Lemma 2, we obtain
$$
m=d(x_b,C_m)=d(x_b,C_m[x_{a+1},x_{b-1}])+d(x_b,C_m[x_{b+1},x_{a-1}])
$$ 
$$
\leq |C_m[x_{a+1},x_{b-1}]|+|C_m[x_{b+1},x_{a-1}]|+2=m.
$$
Therefore
$$
d(x_b,C_m[x_{a+1},x_{b-1}])=|C_m[x_{a+1},x_{b-1}]|+1\\\ \hbox{and} \\\  d(x_b,C_m[x_{b+1},x_{a-1}])=|C_m[x_{b+1},x_{a-1}]|+1.
$$
Since $x_b$ cannot be inserted  into $C_m[x_{b+1},x_{b-1}]$, the last two equalities together with Lemma 2(ii) imply that $x_bx_{a+1}$ and $x_{a-1}x_b\in E(D)$. 

Similarly, one can show that $x_{b-1}x_a$ and $x_ax_{b+1}\in E(D)$. Claim 3 is proved. \fbox \\\\

\noindent\textbf{Claim 4}. {\it Suppose that two distinct     
 vertices, say $x_{m}$ and $x_{q}$, of $\{x_{t_1},x_{t_2},\ldots , x_{t_n}\}$ are not adjacent. Then $\{x_{t_1},x_{t_2},\ldots , x_{t_n}\}$ 
 is an independent
set and for every vertex $y\in \{x_{t_1},x_{t_2},\ldots , x_{t_n}\}$ the following holds
$$
O(x,C_m)=O(y,C_m) \quad   \hbox{and} \quad  I(x,C_m)=I(y,C_m).
$$}
\noindent\textbf {Proof of Claim 4}.  Then $2\leq q\leq m-2$.  By Claim 3 and (18), the arcs
$$
x_{q-1}x, \, xx_{q+1}, \, x_{m-1}x, \, xx_{1}, \, x_{q}x_1,\, x_{m-1}x_{q},\,  x_{q-1}x_{m}\, \hbox{and} \, x_{m}x_{q+1}\, \hbox{are in } \,E(D). \eqno(21)
$$

First we prove the following two statements:

\noindent\textbf{(i)} if $x_jx\in E(D)$ with $j\in [1,m-2]$, then $x_mx_{j+1}\notin E(D)$;

\noindent\textbf{(ii)} if $xx_j\in E(D)$ with $j\in [2,m-1]$, then $x_{j-1}x_m\notin E(D)$.

 \noindent\textbf{Proof of (i) and (ii)}. The proof is by contradiction.  \noindent\textbf{(i)}. Assume the opposite  that $x_jx\in E(D)$ with $j\in [1,m-2]$ and $x_mx_{j+1}\in E(D)$. Now using (21), we obtain that: 
 
if $j\in [1,q-1]$, then $C_{m+1}= x_1x_2\ldots x_jxx_{q+1}\ldots x_mx_{j+1}\ldots x_qx_1$; and

 if $j\in [q+1,m-2]$, then $C_{m+1}= x_1x_2\ldots x_{q-1}x_mx_{j+1}\ldots x_{m-1}x_q\ldots x_jxx_1$. 

 \noindent\textbf{(ii)}.  Assume the opposite   that  $xx_j\in E(D)$ with $j\in [2,m-1]$ and $x_{j-1}x_m\in E(D)$. Again using (21), we obtain that: 

if $j\in [2,q-1]$, then $C_{m+1}= x_1x_2\ldots x_{j-1}x_mx_{q+1}\ldots x_{m-1}xx_j\ldots x_qx_1$; and 
 
 if $j\in [q+1,m-1]$, then $C_{m+1}= x_1x_2\ldots x_{q-1}xx_j\ldots x_{m-1}x_q\ldots x_{j-1}x_mx_1$.

 Thus, in each case we have a contradiction since $C_m$ is a longest cycle in $D$. Therefore, \textbf{(i)} and \textbf{(ii)} both are true.  \fbox \\\\

Now we return to proof of Claim 4.

Let $x_j$ be an arbitrary vertex of $\{x_{t_1},x_{t_2},\ldots , x_{t_n}\}$ other than $x_m$ and $x_q$  (here we assume that $n\geq 3$). Notice that $j\in [2,m-2]$. By (18), $x_{j-1}x$ and 
$xx_{j+1}\in E(D)$. Therefore from (i) and (ii) it follows that $x_mx_j\notin E(D)$ and $x_jx_m\notin E(D)$, 
i.e., $x_m$ and $x_j$ are non-adjacent. Since $x_m$ and $x_j$ are  two arbitrary vertices of $\{x_{t_1},x_{t_2},\ldots , x_{t_n}\}$, we can conclude that any two distinct vertices of  $\{x_{t_1},x_{t_2},\ldots , x_{t_n}\}$ are not adjacent. Thus the first part of Claim 4 is proved, i.e., $\{x_{t_1},x_{t_2},\ldots , x_{t_n}\}$ is an independent set. 

To complete the proof of Claim 4, it remains to show that 
$$
O(x,C_m)=O(y,C_m) \quad   \hbox{and} \quad  I(x,C_m)=I(y,C_m) \quad \hbox{for all} \quad y\in \{x_{t_1},x_{t_2},\ldots , x_{t_n}\}.
$$
Let $x_a$ and $x_b$ be two arbitrary distinct vertices of  $\{x_{t_1},x_{t_2},\ldots , x_{t_n}\}$ such that
$$
\{x_{a+l},x_{a+l+1}, \ldots , x_{b-1}\}\rightarrow x\rightarrow \{x_{a+1},x_{a+2}, \ldots , x_{a+l}\},
$$
in other words, if for some $i\in [1,n]$, $x_a=x_{t_i}$, then  $x_b=x_{t_{i+1}}$.
Choose an arbitrary vertex of $\{x_{t_1},x_{t_2},\ldots , \\ x_{t_n}\}$, say $x_m$ (possibly, $x_m=x_a$).   Then using Lemma 2 and the fact that $d(x_m,C_m)=m$ it is not difficult to show that
$$
d(x_m, C_m[x_{a+1},x_{b-1}])=|C_m[x_{a+1},x_{b-1}]|+1 \quad \hbox{and} \quad x_mx_{a+1}, x_{b-1}x_m\in E(D). \eqno(22)
$$
We claim that $x_m$ and every vertex of $\{x_{a+1},x_{a+2}, \ldots , x_{b-1}\}$ are adjacent. Assume that this is not the case. Let $x_m$ and $x_i$ are not adjacent, where $x_i\in \{x_{a+2},x_{a+3}, \ldots , x_{b-2}\}$. 
If $x_i\in \{x_{a+2},x_{a+3}, \ldots , x_{a+l}\}$, then $xx_i\in E(D)$ and , by statement (ii), $x_{i-1}x_m\notin E(D)$. Now using 
Lemma 2, we obtain 
$$
d(x_m, C_m[x_{a+1},x_{b-1}])=d(x_m, C_m[x_{a+1},x_{i-1}])+d(x_m, C_m[x_{i+1},x_{b-1}])$$
$$
\leq |C_m[x_{a+1},x_{i-1}]|+|C_m[x_{i+1},x_{b-1}]|+1= |C_m[x_{a+1},x_{b-1}]|, 
$$
which contradict the equality of (22). Similarly, using statement (i), one can show that $x_m$ and every vertex of 
$\{x_{a+l+1},x_{a+l+2},$ $ \ldots , x_{b-2}\}$ are adjacent.

Thus the vertex $x_m$ and every vertex of $C_m[x_{a+1},x_{b-1}]$ are adjacent. This together with statement (i) (respectively, statement (ii)) implies  that $\{x_{a+l+1},x_{a+l+2}, \ldots , x_{b-1}\}\rightarrow x_m$ 
and $x_mx_i\notin E(D)$ for all
$x_i\in \{x_{a+l+1},x_{a+l+2},$ $ \ldots , x_{b-1}\}$ (respectively,  
 $x_m\rightarrow \{x_{a+1},x_{a+2}, \\ \ldots , x_{a+l-1}\}$ and $x_ix_m\notin E(G)$ for all $x_i\in \{x_{a+1},x_{a+2},\ldots ,$ $ x_{a+l-1}\}$).  Since
 $$d(x_m, C_m[x_{a+1},x_{b-1}])=|C_m[x_{a+1},x_{b-1}]|+1,$$
 we also have 
$x_m\rightarrow x_{a+l}\rightarrow x_m$. Thus we have proved that $O(x_m,C_m[x_a^+,x_b^-])=O(x,C_m[x_a^+,x_b^-])$ and 
$I(x_m,C_m[x_a^+,x_b^-])=I(x,C_m[x_a^+,x_b^-])$. Therefore, $O(x_m,C_m)=O(x,C_m)$ and $I(x_m,C_m)=I(x,C_m)$.  Since $x_m$ is an arbitrary vertex of $\{x_{t_1},x_{t_2},\ldots , x_{t_n}\}$, we have $O(x,C_m)=O(y,C_m)$ and $I(x,C_m)=I(y,C_m)$ for all $y\in \{x_{t_1},x_{t_2},\ldots , x_{t_n}\}$, which completes the proof of Claim 4.  \fbox \\\\

We  now assume that the fourth statement  of the theorem is not true. Then for some $r\in [2,m-1]$, $D$ contains no cycle of length $r$.  
  Since $d(x,V(C_m))=m$ (recall that $\{x\}=A_1$), it follows that for every $i\in [1,m]$ the following holds
$$
|E(x\rightarrow x_i)|+|E(x_{i+r-2}\rightarrow x)|=1                     \eqno (23)
$$
(for otherwise if $xx_i$ and $x_{i+r-2}x\in E(D)$, then $xx_ix_{i+1}\ldots x_{i+r-2}x$ is a cycle of length $r$, a contradiction).

For any $i\in [1,n]$,  put
$$
F_{2i+1}:=\{x_{t_i+1},x_{t_i+2},\ldots , x_{t_i+l_i-1}\}, \quad 
F_{2i+2}:=\{x_{t_i+l_i+1},x_{t_i+l_i+2},\ldots , x_{t_{i+l}-1}\},
$$
$f_j:=|F_j|$, $y_i:=x_{t_i}$ and $z_{i+1}:=x_{t_i+l_i}$ (all subscripts of $F_j$ are taken modulo $2n$ and all subscripts of $y_i$ and $z_i$ are taken modulo $n$, 
 in particular, $F_{2n+1}=F_1$,  $y_{n+1}=y_1$ and $z_{n+1}=z_1$).\\

\noindent\textbf{Claim 5}. If $|C_m[x_{t_i},x_{t_j}]|=r$, where $i,j\in [1,n]$ and $i\not=j$, then $f_{2i+1}=f_{2j+1}$
and $f_{2i+2}=f_{2j+2}$.

\noindent\textbf {Proof of Claim 5}. We can adjust the notation such that $x_{t_i}=x_{t_1}=x_1$, $x_{t_i+l_i}=x_a$ and $x_{t_2}=x_b$ (possibly, $x_{t_2}=x_{t_j}$). Then $x_{t_j}=x_r$, $F_{2i+1}=F_3=\{x_2,x_3,\ldots , x_{a-1}\}$, $F_{2i+2}=F_4=\{x_{a+1},x_{a+2},\ldots , x_{b-1}\}$ and 
$E(x,\{x_1,x_{b},x_r\})=\emptyset$. In particular, $f_3=a-2$ and $f_4=b-a-1$. 
From the definitions of $F_3$, $F_4$ and (19) it follows that 
$$
\{x_{a}, \ldots , x_{b-1}\}\rightarrow x\rightarrow \{x_2,\ldots , x_a\} \quad  \hbox{and}$$ $$ E(\{x_2,\ldots , x_{a-1}\}\rightarrow x)=
E(x\rightarrow \{x_{a+1},\ldots , x_{b-1}\})= \emptyset.
$$
Using this and  (23), we obtain that  $x_{r-1}x\in E(D)$, 
$$\{x_{a+r-1}, \ldots , x_{b+r-2}\}\rightarrow x \quad \hbox{ and} \quad  
E(\{x_r,x_{r+1},\ldots , x_{a+r-2}\}\rightarrow x)=\emptyset.$$
 The last equality together with Proposition 1(i) implies that 
 $x\rightarrow \{x_{r+1},\ldots , x_{a+r-1}\}$. Since $x_{a+r-1}x\in E(D)$, we have that $|E(x,x_{a+r-1})|=2$. So,
 $ F_{2j+1}=\{x_{r+1},\ldots , x_{a+r-2}\}$. Therefore, $f_{2j+1}=a-2=f_3$.

Since $xx_{b+1}\in E(D)$ (by (18)), from (23) it follows that $x_{b+r-1}x\notin E(D)$. On the other hand, from $x_{b+r-2}x\in E(D)$ and the maximality of $C_m$ it follows that $xx_{b+r-1}\notin E(D)$. 
Therefore, $E(x,x_{b+r-1})=\emptyset$. 
This together with $ \{x_{a+r-1},\ldots , x_{b+r-2}\}\rightarrow x$ and $|E(x,x_{a+r-1})|=2$ implies that 
$F_{2j+2}=\{x_{a+r},\ldots ,$ $ x_{b+r-2}\}$. Hence, $f_{2j+2}=b-a-1=f_4$. This completes the proof of Claim 5.   \fbox \\\\

We are now ready to prove the fourth statement of the theorem.

Without loss of generality, we may assume that $f_1=max\{f_j| 1\leq j\leq 2n\}$. For every $l\in [1,n]$ we will consider the following intervals of integers 
$$
I_1:=[2,f_1+f_2+2] \quad \hbox{and} \quad I_l:=\left[\sum_{i=2}^{2l-1}f_i+2l,\sum_{i=1}^{2l}f_i+2l\right], \quad \hbox{if} \quad l\in [2,n].
$$
Note that 
$$
\sum_{i=2}^{2l-1}f_i+2l=|C_m[z_1,z_l]|+1, \quad \sum_{i=1}^{2l}f_i+2l=|C_m[y_n^+,y_l^-]|+1=|C_m[y_n^+,y_l]|
$$ and  $m=right\{I_n\}$. From the definition of $F_i$ and (23) it follows that $r\notin \cup_{i=1}^nI_i$. It is easy to see that for every $j\in [2,n]$ the following holds
 $$
right\{I_{j-1}\}+1\geq left\{I_j\}-1 \quad \hbox{and} \quad
right\{I_j\}\geq right\{I_{j-1}\}+2,     \eqno (24)
$$
since $f_1$ is maximal.
Now, since  $r\notin \cup_{i=1}^nI_i$, from (24) it follows that there exists an integer $q\in [1,n-1]$ such that 
$$
right\{I_{q}\}+1\leq r\leq left\{I_{q+1}\}-1.      
$$
This together with maximality of $f_1$ gives $f_{1}=f_{2q+1}$ and the last inequalities, in fact, are equalities, i.e.,
 $$
r=right\{I_{q}\}+1=\sum_{j=1}^{2q}f_j+2q+1=left\{I_{q+1}\}-1=\sum^{2q+1}_{j=2}f_j+2q+1.    \eqno (25)
$$ 
Therefore,
$$
r-1=|C_m[y_n^+,y_q]|=|C_m[z_1,z_{q+1}^-]|,\quad \hbox{i.e,} \quad r=|C_m[y_n,y_q]|.  \eqno (26)
$$ 
From (26) and Claim 5 it follows that  $f_1=f_{2q+1}$ and $f_2=f_{2q+2}$. Therefore, $|C_m[y_1,y_{q+1}]|\\=r$. 
Again applying Claim 5, we obtain that $f_3=f_{2q+3}$ and $f_4=f_{2q+4}$. Proceeding in this manner, one can prove that $f_i=f_{i+2q}$ for all $i\in [1,2n]$. This and (25) imply that for every $l\in [1,n]$ the following holds
$$
r=\sum^{2q+2l-2}_{j=2l-1}f_j+2q+1.    \eqno (27)
$$

The remainder of the proof is divided into two cases. \\

\noindent\textbf{Case}. {\it For some $i$ and $j$ ($1\leq i<j\leq n$) the vertices $y_i$ and $y_j$ are not adjacent}. 

Then, by Claim 4,
we have that $\{x,y_1,y_2,\ldots , ,y_n\}$ is an independent set and 
$$
O(x,C_m)=O(y_a,C_m)\quad   \hbox{and} \quad  I(x,C_m)=I(y_{a},C_m) \eqno (28)
$$
for all $a\in [1,n]$. In particular, $y^-_qy_n$ and $ y_nz_n\in E(D)$ since $y^-_qx\in E(D)$ and $xz_n\in E(D)$. Recall that 
$r-1=|C_m[y_n^+,y_q]|$ by (26).
Therefore, if $f_1\geq 1$ (i.e., $xx_2\in E(D)$, we assumed that $y_n=x_m$), then the cycle 
$xC_m[x_2,y^-_q]y_n$ $z_nx$ is a cycle of length $r$, which contradicts the  assumption that $D$ contains no cycle of length $r$. Assume therefore that $f_1=0$. From the maximality of $f_1$ it follows that $f_i=0$ for all $i\in [1,2n]$. This means that
$$
O(x,C_m)= I(x,C_m)=\{x_1,x_3,\ldots , x_{m-1}\}\quad \hbox{and} \quad \{y_1,y_2,\ldots , y_n\}=\{x_2,x_4,\ldots , x_m\}.                \eqno(29)
$$
Since $D$ is strong and $C_m$ is a longest cycle in $D$, from (29) it follows that $A_{2,h}=\emptyset$, i.e., $m=p-1$ and
$V(D)=\{x,x_1,x_2,\ldots ,  x_m\}$ (for otherwise, if $A_{2,h}\not=\emptyset$, then using (29) and $E(A_h\rightarrow C_m)\not= \emptyset$, it is easy to describe a cycle of length greater than $m$). Therefore, $p$  odd and $n=\lfloor p/2\rfloor$.
 Now from (28), (29) and condition $M_0$ it follows that $d(x)=d(x_i)=p-1$ and 
$d(x_i,\{x_j\})=2$ for all $i\in \{x_2,x_4,\ldots ,  x_{m}\}$ and $j\in \{x_1,x_3,\ldots ,  x_{m-1}\}$. 
Now, since 
$\{x,x_2,x_4,\ldots , x_{p-1}\}$ is an independent set and since
$D$ contains no cycle of length $r$, it is not difficult to show that $\{x_1,x_3,\ldots, x_{p-2}\}$ also is an independent set. Therefore, $D\cong K^*_{\lfloor p/2\rfloor,\lfloor p/2\rfloor+1}$. The proof of the fourth statement  of the theorem in this case is complete.\\\\

\noindent\textbf{Case}. {\it  Any pair of vertices $\{y_1,y_2,\ldots , y_n\}$ are adjacent}.\\ 

For this case we need to prove the following claim.\\

 \noindent\textbf{Claim 6. (i)}.  {\it If $y_jy_i\notin E(D)$, where $1\leq i\not= j\leq n$, then $y_iy_j\in E(D)$ and $y_iy^+_j\in E(D)$.}

\textbf{(ii)}. {\it For every integer $l\geq 1$, $y_ny_{lq}\in E(D)$ and $y_ny^+_{lq}\in E(D)$ (recall that 
$r=|C_m[y_n,y_q]|$ by (26)).}\\

\noindent\textbf{ Proof}. \textbf{(i)}. Assume that $y_jy_i\notin E(D)$. Then $y_iy_j\in E(D)$ since the vertices $y_i$
and $y_j$ are adjacent. Since $y^-_ix\in E(D)$, $xy^+_i\in E(D)$ (by (18)) and   $C_m$ is a longest cycle in $D$, it follows that the vertex $y_i$ cannot be inserted into $C_m[y^+_i,y^-_i]$. Now using Lemma 2 and the assumption that $y_jy_i\notin E(D)$, we obtain
$$
d(y_i,C_m[y^+_i,y_j])\leq |C_m[y^+_i,y_j]| \quad \hbox{and} \quad d(y_i,C_m[y^+_j,y^-_i])\leq |C_m[y^+_j,y^-_i]|+1
$$
This together with $d(y_i,C_m)=m$ (by (20)) implies that 
$$d(y_i,C_m[y^+_j,y^-_i])=|C_m[y^+_j,y^-_i]|+1.$$ 
Hence, by Lemma 2 we have, $y_iy^+_j\in E(D)$ since $y_i$ cannot be inserted into $C_m[y_j^+,y_i^-]$.

\textbf{(ii)}. We proceed by induction on $l$. Let $l=1$. If $y_qy_n\in E(D)$, then, since $r=|C_m[y_n,y_q]|$ (by (26), the cycle $y_qy_nC_m[y_n^+,y_q^-]y_q$ has length equal to $r$, which contradicts  our assumption that $D$ contains no cycle of length $r$. Assume therefore that $y_qy_n\notin E(D)$.
Then, by Claim 6(i),   $y_ny_q\in E(D)$ and $y_ny^+_q\in E(D)$. Therefore for $l=1$ Claim 6(ii) is true.

 Assume now that
 $y_ny_{(l-1)q}\in E(D)$ and $y_ny^+_{(l-1)q}\in E(D)$, where $l\geq 2$, and prove that $y_ny_{lq}$ and $y_ny^+_{lq}\in E(D)$.
Assume that $y_{lq}y_{n}\in E(D)$. If in (27) instead of $l$ replace $q(l-1)+1$, then we get 
$$
r-1=\sum^{2ql}_{j=2q(l-1)+1}f_j+2q.
$$ 
This means
  that $|C_m[y^+_{q(l-1)},y_{lq}]|=r-1$, i.e.,
   $|C_m[y_{q(l-1)},y_{lq}]|=r$. 
   Since $y_ny^+_{(l-1)q}\in E(D)$, we have that 
  
if $y^+_{lq}\in C_m[y^+_{(l-1)q},y_n^-]$, then $y_nC_m[y_{(l-1)q}^+,y_{lq}]y_n$ is a cycle of length $r$; and
 
if $y_{lq}\in C_m[y_{n}^+,y_{(l-1)q}^-]$, then $xC_m[y_{n}^+,y_{lq}]y_nC_m[y_{(l-1)q}^+,y_{n}^-]x$ is a cycle of length $r$. 
In both cases we have a contradiction. 

We may therefore  assume that  $y_{lq}y_{n}\notin E(D)$. Then from Claim 6(i) it follows that $y_ny_{lq}\in E(D)$ and $y_ny^+_{lq}\in E(D)$. This completes the proof of Claim 6.
     \fbox \\\\

 It is not difficult to see that there exists an integer $l\geq 2$ such that $y_n=y_{lq}$. To see this it suffice to notice that $y_{nq}=y_n$ since the subscripts  of $y_i$ are considered modulo $n$. Therefore, by Claim 6(ii) we have that
$y_ny_{(l-1)q}\in E(D)$ and hence, the cycle $y_nC_m[y_{(l-1)q},y_n^-]y_n$ has length equal to $r$, which is a contradiction. Thus  the fourth statement of the theorem also is proved. This completes the proof of the theorem. \fbox \\\\

Note that the Thomassen theorem is a consequence of Theorem D. From Theorem D also it follows the following corollary.\\
 
\noindent\textbf{Corollary}.  {\it Let $D$ be a strongly connected non-Hamiltonian  digraph of order $p$ and with minimum degree at least 
$p-1$ and let $m$ be the length of a longest cycle in $D$. Then} 

(i) {\it $D$ contains a cycle of length $r$ for all $r\in [2,m]$ unless when $p$ is odd  and $D$ is isomorphic to the complete bipartite digraph $ K^*_{\lfloor p/2\rfloor,\lfloor p/2\rfloor+1}$}. 

(ii) {\it If $D$ is strongly 2-connected, then the induced subdigraph $D\langle V(D)\setminus V(C_m)\rangle$ is a transitive tournament, where $C_m$ is an arbitrary longest cycle in $D$.}\\

\noindent\textbf{Note added in proof}. In view of  fourth statement of  Theorem D it is natural to pose the following problems.\\

\noindent\textbf{Problem 3.} Let $D$ be a Hamiltonian digraph satisfying condition $M_0$. Whether $D$ is pancyclic?\\

\noindent\textbf{Problem 4.} Find sufficient conditions    for Hamiltonian digraphs with condition $M_0$ to be pancyclic.

\end{document}